\documentclass[12pt]{article}

\title{The Kauffman skein algebra of a surface at $\sqrt{-1}$}
\usepackage[a4paper]{geometry}
\usepackage[applemac]{inputenc}
\input xy
\xyoption{all}
\usepackage{amsmath}
\usepackage{amsfonts}
\usepackage{amssymb}
\usepackage{stmaryrd}
\usepackage{graphicx,pstricks}
\usepackage{bbm,color}
\usepackage{amsthm}
\DeclareMathOperator{\hol}{Hol}
\DeclareMathOperator{\re}{Re}
\DeclareMathOperator{\acos}{Arcos}

\DeclareMathOperator{\tr}{Tr}
\DeclareMathOperator{\en}{End}

\begin{document}
\author{Julien March\'e}
\date{}
\newcommand{\C}{\mathbb{C}}
\newcommand{\Z}{\mathbb{Z}}
\renewcommand{\H}{\mathbb{H}}
\renewcommand{\hom}{\rm Hom}
\newcommand{\R}{\mathbb{R}}
\newcommand{\N}{\mathbb{N}}
\renewcommand{\lg}{\langle}
\newcommand{\ld}{\rangle}
\newcommand{\boA}{\mathcal{A}}
\newcommand{\boF}{\mathcal{F}}
\newcommand{\boM}{\mathcal{M}}
\newcommand{\boL}{\mathcal{L}}
\newcommand{\bog}{\mathcal{G}}
\newcommand{\boO}{\mathcal{O}}
\newcommand{\boB}{\mathcal{B}}
\newcommand{\boR}{\mathcal{R}}
\newcommand{\boH}{\mathcal{H}}
\newcommand{\boT}{\mathcal{T}}
\newcommand{\su}{\rm{SU}(2)}
\newcommand{\so}{\rm{SO}(3)}
\newcommand{\SL}{\rm{SL}(2,\mathbb{C})}
\renewcommand{\d}{{\rm d}}
\renewcommand{\phi}{{\varphi}}
\renewcommand{\hat}{\widehat}
\newcommand{\ba}[1]{\overline{#1}}
\newcommand{\mat}[1]{\begin{bmatrix} #1 \end{bmatrix}}
\maketitle
\newtheorem{thm}{Theorem}[section]
\newtheorem{defn}[thm]{Definition}
\newtheorem{lem}[thm]{Lemma}
\newtheorem{quest}[thm]{Question}
\newtheorem{prop}[thm]{Proposition}
\newtheorem{rem}[thm]{Remark}
\newtheorem*{thm*}{Theorem}
\newtheorem*{prop*}{Proposition}

\abstract{We study the structure of the Kauffman algebra of a surface with parameter equal to $\sqrt{-1}$. We obtain an interpretation of this algebra as an algebra of parallel transport operators acting on sections of a line bundle over the moduli space of flat connections in a trivial \su-bundle over the surface. We analyse the asymptotics of traces of curve-operators in TQFT in non standard regimes where the root of unity parametrizing the TQFT accumulates to a root of unity. We interpret the case of $\sqrt{-1}$ in terms of parallel transport operators.
}
\section{Introduction}
The Kauffman skein algebra of a surface is a fascinating object: its very simple and combinatorial definition hides deep and unobvious links with representation varieties of the surface in subgroups of SL(2,$\C$), their symplectic geometry, deformation and exact quantizations. Let us recall briefly some of these relations starting from the definition (see \cite{bfk} for a nice introduction).
Let $\Sigma$ be an oriented compact surface without boundary and $A$ be a formal parameter. 
We set $K(\Sigma,A)$ to be the free $\Z[A,A^{-1}]$-module generated by isotopy classes of banded links in $\Sigma\times [0,1]$ modulo the following local relations:
\begin{figure}[width=8cm,height=3cm]
\begin{center}
\begin{pspicture}(-2,0)(3,3)
\includegraphics{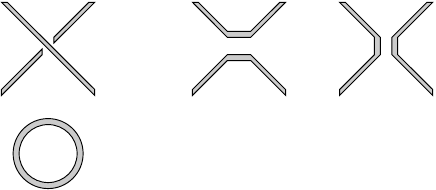}
\put(-5.2,2.3){$=A$}
\rput(-1.9,2.5){$+A^{-1}$}
\rput(-4,0.7){$=(-A^2-A^{-2})\quad \emptyset$}
\end{pspicture}
\caption{Kauffman relations}\label{kauf}
\end{center}
\end{figure}

The stacking of two banded links define on $K(\Sigma,A)$ a structure of algebra. We define a {\it multicurve} as a 1-submanifold of $\Sigma$ without component bounding a disc. It is well known that $K(\Sigma,A)$ is a free $\Z[A,A^{-1}]$-module generated by isotopy classes of multicurves (see \cite{ps}). For any $u\in\C^*$, we set $K(\Sigma,u)=K(\Sigma,A)\otimes_{\Z[A,A^{-1}]} \C$ where $A$ acts on $\C$ by $u$. All these algebras are naturally isomorphic as vector spaces when we identify the basis of multicurves.

It is shown in \cite{bull,ps} that $K(\Sigma,-1)$ is isomorphic to the algebra of regular functions on the representation variety $\boM(\Sigma,\SL)$. This isomorphism sends a simple curve $\gamma$ to the function $f_{\gamma}$ defined by $f_{\gamma}(\rho)=-\tr\rho(\gamma)$. 

Moreover, $K(\Sigma,-e^h)$ gives an algebra structure on $K(\Sigma,-1)[[h]]$ which deforms the preceding one and may be interpreted as a star-product whose first order is given by the Poisson bracket of the well-known symplectic structure on $\boM(\Sigma,\su)$ (see \cite{goldman2}). Finally, if $A$ is a primitive $2p$-th root of unity, the algebra $K(\Sigma,A)$ has a quotient which may be viewed as the space of endomorphisms of $V_p(\Sigma)$, the quantization of level $p/2-2$ of the moduli space $\boM(\Sigma,\su)$ defined in \cite{BHMV}.

All these instances of the Kauffman skein algebra should be the manifestation of a geometric meaning of the full structure, or at least of its specializations at roots of unity. 
This is not the case at the moment of our understanding and what we propose to study is yet another manifestation of these expected relations.

Consider the Kauffman skein algebra $K(\Sigma,-i)$. This is a non-commutative algebra for which we propose the following interpretation: 
let $\boM(\Sigma,\su)$ be the moduli space of flat connections on the trivial principal bundle $\Sigma\times \su\to\Sigma$. This space has a natural integral symplectic form on its smooth part and a natural prequantum bundle $\boL\to\boM(\Sigma,\su)$ with hermitian structure and connection.
Given a simple curve $\gamma$, let $F_{\gamma}$ be the function defined on $\boM(\Sigma,\su)$ by $F_{\gamma}(\rho)=\frac{1}{\pi}\acos\frac{1}{2}\tr \rho(\gamma)$.
Where it is smooth, this function generate a 1-periodic hamitonian flow. Thanks to the connection on $\boL$, the flow extends to the line bundle: let us denote it by $\Phi_{\gamma}^{t}$.

We define an operator $\boO_{\gamma}$ acting on $L^2(\boM(\Sigma,\su),\boL)$ by the formula $(\boO_{\gamma}s)(\rho)=\Phi_{\gamma}^{1/2} s(\Phi_{\gamma}^{-1/2}\rho)+\Phi_{\gamma}^{-1/2} s(\Phi_{\gamma}^{1/2}\rho)$. In unformal terms, $\boO_{\gamma}s$ is the parallel transport of $s$ along the hamiltonian flow of $F_{\gamma}$ during a time $1/2$ plus the parallel transport of the same flow during a time $-1/2$. 
Our main result is the following:
\begin{thm*}
The map $K(\Sigma,-i)\to {\rm End\,} L^2(\boM(\Sigma,\su),\boL)$ is an injective morphism of algebras.
\end{thm*}
As an intermediate step, we use a description of the skein algebra at 4-th roots of unity very similar to the one obtained in \cite{si}.
Moreover, there is a natural and geometric trace on the operator algebra generated by all $\boO_{\gamma}$: this will induce a trace on $K(\Sigma,-i)$. By computing the limit of traces of curve operators in TQFT for asymptotic regimes converging to any roots of unity, we identify our  geometric trace with a limit of TQFT traces. In this way, we identify the Fell limit of TQFT whose defining primitive root converges to $-i$ with a geometric representation of the mapping class group. This part is very similar to the article \cite{mn} and so will not be developed in full length. Let us comment the structure of the article.

\begin{itemize}
\item[-]
The first part of the article is a discussion about the case when $\Sigma$ is a torus. In that case, everything is known, the structure of the Kauffman bracket is understood and the moduli space and its prequantum bundle are particularly simple: we give the full description of both and a proof of the proposition in that case. The reason is first to give to the reader an easy picture of what happens in general, and second to explain the way we arrived to the result explained  in the article. 
\item[-]
In the second part, we identify the Kauffman skein algebra at $-i$ with an algebra which is a mixture of the Kauffman algebra at -1 with an Heisenberg-type algebra.
\item[-] 
Next, we give many details about moduli spaces of flat connections to state our geometric interpretation. In particular, we define the Heisenberg group via the lifting of the gauge group with values in $\so$, and give explicit formulas for Hamiltonian flows and their lift to the prequantum bundle. No originality is claimed here: we adapted standard results to our purpose.
\item[-]
Finally, we compute some limit of traces in TQFT and identify the limit when the roots of unity involved in the limit converge to $-i$.
\end{itemize}

{\bf Acknowledgements}
I would like to thank Laurent Charles for discussing much of the matter of this article, Alex Bene for the proof of the lemma 3.3 and Gregor Masbaum for his interest and remarks.

\section{The torus case}
Let $\Sigma$ be an oriented torus with a choice of oriented meridian $m$ and parallel $l$ such that $m\cdot l=1$.

A basis of $K(\Sigma,A)$ is obtained from simple curves and their powers. We denote by $(p,q)$ the curve $pm+ql$ where $p$ and $q$ are relatively prime integers.

Following \cite{torus} and ideas of \cite{sal}, let us denote by $\boT$ the non-commutative algebra defined by 
$\boT=\boR \lg L^{\pm 1},M^{\pm 1}\ld/(LM=A^2ML)$. Let $\boT^{\sigma}$ be the subalgebra of $\boT$ invariant under the involution $\sigma$ defined by $\sigma(M^pL^q)=M^{-p}L^{-q}$. 
Then, Frohman and Gelca showed that the map $\Phi:K(\Sigma,A)\to \boT^{\sigma}$ defined by $\Phi(p,q)=A^{pq}(M^pL^q+M^{-p}L^{-q})$ is an isomorphism of algebra. We can interpret this theorem in a geometric way by using the representation variety of $\Sigma$ in $\su$.

More precisely, given $\rho:\pi_1(\Sigma)\to \su$ we can find a basis where $\rho(m)$ and $\rho(l)$ are simultaneously diagonalized as these matrices commute.
That is, we can find $\alpha, \beta \in \R$ such that up to conjugation, one has 
$$\rho(m)=\mat{e^{i\pi \alpha}&0\\0&e^{-i\pi\alpha}}\text{ and }\rho(l)=\mat{e^{i\pi \beta}&0\\0&e^{-i\pi\beta}}.$$
The parameters $\alpha$ and $\beta$ are well-defined only modulo 2 and the representation $\rho$ does not change if we replace $(\alpha,\beta)$ with $(-\alpha,-\beta)$. Hence, we can identify $\boM(\Sigma,\su)$ with $(\R^2/2\Z^2)/(x\sim -x)$. The symplectic form $\omega$ corresponds to $\frac{1}{2}\d \alpha\wedge \d \beta$ such that the volume of $\boM(\Sigma,\su)$ is equal to 1.

Next, for a simple curve $\gamma$ in $\Sigma$, we consider the function $F_{\gamma}$ defined on $\boM(\Sigma,\su)$ by the formula $F_{\gamma}(\rho)=\frac{1}{\pi}\acos (\frac{1}{2}\tr \rho(\gamma))$. We have for instance $F_{m}(\alpha,\beta)=|\alpha| \mod 1$ and $F_{l}(\alpha,\beta)=|\beta|\mod 1$. These functions are differentiable in a dense open set of $\boM(\Sigma,\su)$ and their hamiltonian flows noted  $\Phi^t_{\gamma}$ are 1-periodic. One has for instance $\Phi^t_m(\alpha,\beta)=(\alpha,\beta+2t)$ and $\Phi^t_l(\alpha,\beta)=(\alpha-2t,\beta)$.

Let us fix a pre-quantum bundle $\boL$ over $\boM(\Sigma,\su)$ that is a complex line bundle with hermitian structure and connexion whose curvature is the symplectic form. 
We can take a quotient of the trivial bundle $\R^2\times\C$ by the following action of $\Z^2\rtimes \Z/2$. The $\Z^2$ subgroup acts by $(m,n).(\alpha,\beta,z)=(\alpha+2m,\beta+2n,\exp(i\pi(\alpha n-\beta m))z)$ and the generator of $\Z/2$ transforms $(\alpha,\beta,z)$ into $(-\alpha,-\beta,z)$.
The 1-form $\frac{2i\pi}{4}(\alpha\d\beta-\beta\d\alpha)$ gives a connection on that bundle with curvature $\omega$. Thanks to that connection, we can extend the flow $\Phi_{\gamma}^t$ in $\boM(\Sigma,\su)$ to a 1-parameter flow in $\boL$: we denote it also by $\Phi_{\gamma}^t$.
Hence, we construct an operator $\Psi_{\gamma}^t$ that acts on sections of $\boL$ via the formula $(\Psi_{\gamma}^t s)(\rho)=\Phi_{\gamma}^t s(\Phi_{\gamma}^{-t}\rho)$. This operator is not defined when $F_{\gamma}$ is 0 or 1, nevertheless we can make sense of these operators as acting in $L^2(\boM(\Sigma,\su),\boL)$.

If we do parallel transport along a null-homotopic loop $\gamma$ in $\boM(\Sigma,\su)$, then the holonomy in the prequantum bundle is $\exp(2i\pi\int_D\omega)$ where $\omega$ is the symplectic form and $D$ is a disc bounding $\gamma$. 
This standard fact implies that we have $\Psi^t_l\Psi^t_m=e^{-4i\pi t^2}\Psi_m^t\Psi_l^t$ as the operator $\Psi^t_m\Psi^t_l\Psi^{-t}_m\Psi^{-t}_l$ is the parallel transport around a square of symplectic area $2t^2$.
Moreover, for $\gamma$ a curve with slope $(p,q)$, one has $\Psi^t_{\gamma}=e^{-2i\pi t^2 pq}(\Psi_m^t)^p(\Psi_l^t)^q$.
These two formulas show that we may interpret the algebraic structure of $K(\Sigma,A)$ as the composition of parallel transport operators in such a way that $\gamma$ corresponds to the operator $\Psi_{\gamma}^t+\Psi_{\gamma}^{-t}$. The value of $t$ is linked with $A$ with the relation $A=\exp(-2i\pi t^2)$. This is the reason why we will consider $t=1/2$ and $A=\exp(-i\pi/2)$ in what follows.

\section{An algebra isomorphic to the Kauffmann skein algebra at $\sqrt{-1}$}

\begin{defn}
Let $\boA$ be the free $\C$-algebra generated by the symbols $[\gamma]$ for $\gamma\in H_1(\Sigma,\Z)$ with the  following relations: for any $\gamma,\delta \in H_1(\Sigma,\Z)$ one has 
\begin{equation}\label{multiplication}
[\gamma]^2=1\quad {\rm and} \quad[\gamma][\delta]=i^{-\gamma\cdot\delta}[\gamma+\delta]
\end{equation}
\end{defn}
This algebra is graded by the abelian group $H^1(\Sigma,\Z_2)$ where the grading of $[\gamma]$ is the class of $\gamma$ mod 2. Moreover, the dimension of each homogeneous piece is 1, generated by any integral cohomology class projecting onto the corresponding mod 2 class.
 We remark that for any value of $u\in\C^*$, the Kauffman skein algebra $K(\Sigma,u)$ is   also graded by $H^1(\Sigma,\Z_2)$ as the Kauffman relations preserve the mod 2 homology class of the links.
Hence, we can consider the tensor product $K(\Sigma,-1)\otimes \boA$ as a graded algebra. We denote by $\boA(\Sigma)$ the sub-algebra of elements with grading 0. Informally, we twist the algebra $K(\Sigma,-1)$ by tensoring any element $\gamma$ with its class $[\gamma]$ in $\boA$. 
The aim of this section is to show the following result:

\begin{thm}\label{combi}
The map $\phi:K(\Sigma,-i)\to \boA(\Sigma)$ defined by $\phi(\gamma)=(-1)^{n(\gamma)}\gamma\otimes [\gamma]$ is an isomorphism of algebras where $\gamma$ is a multicurve on $\Sigma$ and $n(\gamma)$ is the number of its components.
\end{thm}
Notice that using the formulas \eqref{multiplication}, the element $[\gamma]$ is defined for any multicurve and does not depend on an orientation of $\gamma$.
This proposition is an important step of the geometric interpretation: it is the only place where Kauffman relations are dealt with. The  proof will rely on the following lemma.

\begin{figure}
\begin{center}
\begin{pspicture}(-2,-0.5)(2,1)
\includegraphics{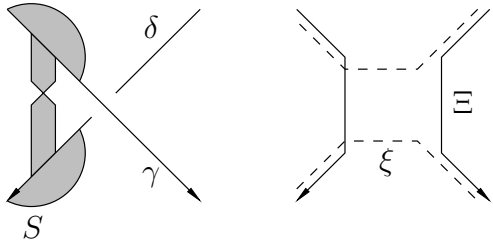}
\put(-6.2,-0.4){$S$}
\rput(-4.5,0.4){$\gamma$}
\rput(-4.5,2.4){$\delta$}
\rput(-1.4,0.6){$\xi$}
\rput(-0.4,1.4){$\Xi$}
\end{pspicture}
\caption{Non-seifert crossings and Moebius bands}\label{moebius}
\end{center}
\end{figure}

\begin{lem}\label{alex}
Let $G$ be a finite graph such that in the neighborhood of any vertex, the edges incoming to that vertex have a cyclic order. We decompose the edges of $G$ in two parts: $E(G)=E_h\coprod E_m$. The edges of $E_h$ will be said to be of type handle whereas the edges of $E_m$ will be of type moebius.

We construct a surface $S$ from these data in the following way: take a family of oriented discs parametrized by vertices of $G$. For all edges, we attach a rectangle to the corresponding discs such that the cyclic orientation of the vertices is respected. The rectangle should respect the orientations of the discs if the edge has type handle and should not respect them if the edge has type Moebius. In the sequel, we will refer to rectangles of the second type as Moebius bands.

Orient the boundary of $S$ in an arbitrary way. Let $n$ be the number of boundary components and $m$ be the number of Moebius bands whose sides are oriented in the same direction. Then the following formula holds:
$$n+m+\chi(S)=0\mod 2.$$
\end{lem}
\begin{proof}
We replace Moebius bands with parallel sides with bands of type handle and cut the Moebius bands with opposite orientations of the sides. In the first case, the Euler characteristic does not change but the number of boundary components differs by one, in the second case, the Euler characteristic and the number of components differ by one. 
Hence, during this operation, the parity of the quantity $n+m+\chi(S)$ did not change but $S$ became an orientable surface. We are left with the equality $n+\chi(S)=0\mod 2$ for $S$ an orientable surface which is obvious.

Let us prove the theorem. It is clear that $\phi$ is an isomorphism of vector spaces because it sends the basis consisting of isotopy classes of multicurves to the same basis of $\boA(\Sigma)$. We also remark that the trivial curve is sent to the element -2 in $\boA(\Sigma)$ which agrees with the first Kauffman relation.

Let $\gamma$ and $\delta$ be two multicurves and suppose that they intersect transversally. Let $\xi$ be any smoothing of $\gamma\cup\delta$: we denote by $c(\xi)$ the number of positive resolutions minus the number of negative resolutions when considering that $\gamma$ lies above $\delta$. Proving the theorem is equivalent to showing that the following formula holds:
\begin{equation}\label{eq1}
\phi(\gamma)\phi(\delta)=\sum\limits_{\xi} i^{-c(\xi)}\phi(\xi).
\end{equation}
Suppose that $\gamma$ and $\delta$ are oriented and denote by $\Xi$ the Seifert smoothing of $\gamma\cup\delta$. As the homology class of $\Xi$ represents the sum of the integer homology classes of $\gamma$ and $\delta$, we have $[\gamma][\delta]=i^{-\gamma\cdot\delta}[\Xi]$. Moreover, we know that $\xi$ and $\Xi$ have the same homology class mod 2 and so the classes $[\xi]$ and $[\Xi]$ are proportional in $\boA$. By the multiplication formula \eqref{multiplication}, we obtain $[\xi]=i^{\Xi\cdot\xi}[\Xi]$ (notice that $\Xi\cdot\xi$ is even).

Let us identify a multicurve with its image in $K(\Sigma,-1)$: the formula \ref{eq1} becomes 

\begin{equation}\label{eq2}
 i^{-\gamma\cdot \delta}(-1)^{n(\gamma)+n(\delta)}\gamma\delta\otimes[\Xi]=\sum\limits_{\xi}  i^{-c(\xi)+\Xi\cdot\xi} (-1)^{n(\xi)}\xi\otimes[\Xi].
\end{equation}
 But in $K(\Sigma,-1)$ one already has 
 $$\gamma\cdot\delta=\sum\limits_{\xi}(-1)^{c(\xi)}\xi.$$
 Hence, to prove the theorem, it is sufficient to show that for all smoothings $\xi$ one has $$i^{-\gamma\cdot \delta}(-1)^{n(\gamma)+n(\delta)}=i^{c(\xi)+\Xi\cdot\xi} (-1)^{n(\xi)}.$$ 
The quantity $\gamma\cdot\delta$ is a sum of 1 and -1 over the set of crossings of $\gamma$ and $\delta$, just as $c(\xi)$. We remark that these signs are equal when the smoothing is Seifert and are opposite when it is not, thus one has $\gamma\cdot\delta-c(\xi)=2 ns(\xi)$ if we denote by $ns(\xi)$ the number of non-Seifert smoothings of $\xi$. Knowing that the number of crossings mod 2 is equal to $\gamma\cdot\delta$ one has to prove the following formula:

\begin{equation}
ns(\xi)+\frac{1}{2}\Xi\cdot \xi +n(\xi)+n(\gamma)+n(\delta)+\gamma\cdot\delta=0\mod 2.
\end{equation}
The only number which remain to be computed is $\Xi\cdot\xi$. We use the fact that $\Xi$ is globally oriented as a Seifert smoothing. Far away from the crossing points we isotope $\xi$ to a curve parallel to $\Xi$ and slightly on its left with respect to the orientation of $\Xi$ as in the right hand side of figure \ref{moebius}. We notice that $\Xi$ and $\xi$ may cross only in the neighborhood of non-seifert crossings and that in these neighborhoods, the curves cross in two points. To compute the intersection, we need to orient $\xi$ and count crossings where $\xi$ and $\Xi$ intersect in two points with the same sign. 

We can interpret the situation in the setting of the lemma \ref{alex} by constructing a surface $S$ as follows. Consider for each component of $\gamma$ and $\delta$ an abstract oriented disc which bounds them. For each intersection point of $\gamma$ and $\delta$ we glue a band between the corresponding discs. The band respects the orientation of the discs if the smoothing of $\xi$ is Seifert at this crossing and reverse it if the smoothing is non-Seifert. The boundary of $S$ is precisely $\xi$ and we remark in figure \ref{moebius} that the boundary curves of a Moebius band are oriented in the same direction if and only if $\xi$ intersects $\Xi$ in two opposite ways in the neighborhood of the crossing. Using the notation of the lemma,  we conclude that $m$ is equal to $ns(\xi)+\frac{1}{2}\Xi\cdot\xi\mod 2$. On the other hand, applying the lemma we get the formula
$$m=n(\xi)+\chi(S)=n(\xi)+n(\gamma)+n(\delta)+\gamma\cdot\delta \mod 2.$$
This last equality ends the proof.
\end{proof}

\section{A geometric interpretation}

This section contains standard material, mainly inspired from \cite{jw1} and \cite{jw2}. The appearance of the Heisenberg group was noticed in \cite{BHMV} and \cite{am}: we introduce it in a different way, very standard in its spirit but for which we could not find explicit reference, so we decided to give all the details. The geometric interpretation of $K(\Sigma,-i)$ is given in the last subsection.

\subsection{The SU(2)-moduli space and its prequantum line bundle}\label{generalite}

Let $\Sigma$ be a connected, compact and oriented surface without boundary. 

We will identify $\su$ with the unit sphere of the set of quaternions $\H$. Hence, the Lie algebra of $\su$ that we denote by $\bog$ is identified with the vector space with basis $i,j$ and $k$. We denote by $(\cdot,\cdot)$ the standard scalar product on $\bog$ such that $i,j,k$ is an orthonormal basis. We set $\lg\cdot,\cdot \ld=\frac{1}{2\pi^2}(\cdot,\cdot)$ as a second scalar product; we will use both of them to avoid constants in the formulas.

We identify the space of connections on the principal bundle $\Sigma\times\su\to\Sigma$ with $\Omega^1(\Sigma,\bog)$ with the convention that $\su$ acts on the right. We define the gauge group $\Gamma_{\su}$ as the group of differentiable functions from $\Sigma$ to $\su$. This group acts on $\Omega^1(\Sigma,\bog)$ via the formula $a^g=g^{-1}ag+g^{-1}\d g$.

The curvature of a connection $a$ is the 2-form $\d a+\frac{1}{2}[a\wedge a]$. We denote by $\Omega^1_{\rm flat}(\Sigma,\bog)$ the set of connections whose curvature is 0. We recall that, fixing a base point on $\Sigma$, the holonomy of based loops gives the identification:

$$\Omega^1_{\rm flat}(\Sigma,\bog)/\Gamma_{\su}\simeq \hom(\pi_1(\Sigma),\su)/\su=\boM(\Sigma,\su).$$

Consider the symplectic form on $\Omega^1(\Sigma,\bog)$ defined by $\omega(a,b)=\int_{\Sigma} \lg a\wedge b\ld$. Following \cite{atiyah}, the action of the gauge group on $\Omega^1(\Sigma,\bog)$ is hamiltonian with moment map given by the curvature. The moduli space $\boM(\Sigma,\su)$ inherits a natural symplectic structure on its smooth part as a symplectic reduction. 

We give a sketch of construction of the prequantum bundle. We refer to \cite{freed} for details. Let $\theta$ be the left-invariant Maurer-Cartan 1-form on $\su$ and $\chi$ be the Cartan 3-form defined by $\chi=\frac{1}{12}\lg\theta\wedge[\theta\wedge \theta]\ld$. This 3-form is left and right invariant and our choice of normalisation is such that $\int_{\su}\chi=1$.

For $g\in\Gamma_{\su}$, we define $W_{\su}(g)\in\R/\Z$ in the following way: let $M$ be an oriented 3-manifold with boundary $\Sigma$ and $\tilde{g}$ be an extension of $g$ to $M$. Such an extension always exists as one has $\pi_1(\su)=\pi_2(\su)=0$. We set
 $$W_{\su}(g)=\int_M \tilde{g}^*\chi$$ This is a well-defined element of $\R/\Z$. In fact, if $N$ is another 3-manifold and $\tilde{g}$ an extension of $g$ to $N$, one has $\int_M\tilde{g}^*\chi-\int_N \tilde{g}^*\chi=\int_{M\cup(-N)}\tilde{g}^*\chi=\deg \tilde{g}^*\in \Z$.

We identify $S^1$ with $\R/\Z$ and set for $a\in \Omega^1(\Sigma,\bog)$ and $g\in\Gamma_{\su}$ 
$$c(a,g)=\frac{1}{2}\int_{\Sigma}\lg g^{-1}ag\wedge g^{-1}\d g\ld -W_{\su}(g).$$

This expression satisfies the cocycle relation $c(a,gh)=c(a,g)+c(a^g,h)$. We give some elements of the proof as we will come back to it in the next section: this relation comes from the formula $\int_{M}(\tilde{g}\tilde{h})^*\chi=\int_M\tilde{g}^*\chi+\int_M\tilde{h}^*\chi-\frac{1}{2}\int_{\Sigma}\lg g^{-1}\d g\wedge \d h h^{-1}\ld$ where $M$ bounds $\Sigma$ and $\tilde{g},\tilde{h}$ are extensions of $g$ and $h$ respectively. This formula is itself the integration of the formula $\mu^* \chi= \pi_1^*\chi+\pi_2^*\chi-\frac{1}{2}\lg\pi_1^*\theta,\pi_2^*\overline{\theta}\ld$ where $\mu$ is the product map from $\su\times\su$ to $\su$, $\pi_1,\pi_2$ are the two projections and $\overline{\theta}$ is the right invariant Maurer-Cartan form.

Now we can define the $S^1$-bundle $\boL$ (which is equivalent to defining an hermitian line bundle) in the following way: we define $\boL$ as the quotient of $\Omega^1_{\rm flat}(\Sigma,\bog)\times S^1$ by the following action of $\Gamma_{\su}$: $(a,\theta)^g=(a^g,\theta+c(a,g))$. The connection is the quotient of the 1-form $\lambda$ given by $\lambda_{(a,\theta)}(b,\eta)=\frac{1}{2}\int_{\Sigma}\lg a\wedge b\ld +\eta$. One can check easily that this form is a connection form equivariant with respect to the action of $\Gamma_{\su}$ on $\Omega^1_{\rm flat}(\Sigma,\bog)\times S^1$.

\subsection{The action of the Heisenberg group}\label{heisenberg}

Let $\Gamma_{\so}$ be the group of differentiable functions from $\Sigma$ to $\so$. The connected component of the identity is isomorphic to $\Gamma_{\su}/\{\pm 1\}$ and the group of connected components is isomorphic to $H^1(\Sigma,\Z_2)$. 
The group $\Gamma_{\so}$ acts on $\Omega^1_{\rm flat}(\Sigma,\bog)$ with the usual formula $a^g=g^{-1}ag+g^{-1}dg$. The quotient is isomorphic to $\boM(\Sigma,\so)$. In this way, we obtain that $\boM(\Sigma,\so)$ is a quotient of $\boM(\Sigma,\su)$ by the action of $H^1(\Sigma,\Z_2)$. We can understand this action in terms of representation spaces in the following way:
an element $\lambda \in H^1(\Sigma,\Z_2)$ may be viewed as a group homomorphism from $\pi_1(\Sigma)$ to $\Z_2\simeq\{\pm 1\}\subset \su$. For $\lambda:\pi_1(\Sigma)\to\{\pm 1\}$ and $\rho:\pi_1(\Sigma)\to\su$ two representations, the representation $\lambda.\rho$ is nothing but the product of the corresponding group homomorphisms from $\pi_1(\Sigma)$ to $\su$.

When trying to extend the action of $H^1(\Sigma,\Z_2)$ on $\boM(\Sigma,\su)$ to the prequantum bundle, we encounter an obstruction which is solved by considering a non trivial extension of $H^1(\Sigma,\Z_2)$ by $\Z_4$, namely the (an) Heisenberg group.

\subsubsection{The extended Wess-Zumino-Witten functional}

Let $g$ be an element of $\Gamma_{\so}$. To define $W_{\so}(g)$ we have to face two difficulties: the first one is that given a 3-manifold $M$ bounding $\Sigma$ it is possible that $g$ does not extend to $M$. The cohomology class of $g$ in $H^1(\Sigma,\Z_2)$ has to be the restriction of a cohomology class of $M$. Nevertheless, it is clear that we can find $M$ such that it is the case (for instance a well chosen handlebody). If we denote by $\tilde{g}$ the extension of $g$ to $M$, we set $W_{\so}(g)=\int_M\tilde{g}^* \chi$. 
The second difficulty is that this number is only defined mod $\frac{1}{2}\Z$. Indeed, given another 3-manifold $N$ such that $g$ extends to $N$ as a map $\tilde{g}$ one has 
$\int_M\tilde{g}^* \chi-\int_N\tilde{g}^* \chi=\int_{M\cup (-N)}\tilde{g}^* \chi=(\deg\tilde{g})\int_{\so}\chi=\frac{1}{2}\deg\tilde{g}$.

The following lemma will be crucial to define the extension of $\Gamma_{\so}$:
\begin{lem}\label{lemme1}
Let $g$ and $h$ be two elements of $\Gamma_{\so}$, then one has 
$$W_{\so}(gh)=W_{\so}(g)+W_{\so}(h)-\frac{1}{2}\int_{\Sigma}\lg g^{-1}\d g \wedge \d hh^{-1}\ld \mod\frac{1}{4}\Z$$
\end{lem}
\begin{proof}
Let $(g_*,h_*):\pi_1(\Sigma)\to\Z_2\times\Z_2$ be the morphism induced by $g$ and $h$ on fundamental groups. It induces a covering of $\Sigma$ with four sheets that we denote by $\check{\Sigma}$. By construction, there are two functions $\check{g}$ and $\check{h}$ such that the following diagram is commutative:
$$\xymatrix{\check{\Sigma}\ar[r]^{\check{g},\check{h}}\ar[d] & \su\ar[d] \\ \Sigma\ar[r]^{g,h} & \so}$$

As $\check{g}$ and $\check{h}$ are two functions with values in $\su$, we can apply the argument of the preceding section to obtain the formula: 
\begin{equation}\label{produit}
W_{\su}(\check{g}\check{h})=W_{\su}(\check{g})+W_{\su}(\check{h})-\frac{1}{2}\int_{\check{\Sigma}}\lg \check{g}^{-1}d\check{g}\wedge \d\check{h}\check{h}^{-1}\ld) \mod \Z
\end{equation}

We can relate $W_{\so}(g)$ and $W_{\su}(\check{g})$ in the following way: let $M$ be a 3-manifold bounding $\Sigma$ on which $g$ extends to a map $\tilde{g}$. 
Then there is a covering $\check{M}$ of $M$ which bounds $\check{\Sigma}$. We obtain a map $\tilde{\check{g}}$ on $\check{M}$ by composing the projection and $\tilde{g}$. Using this extension to compute $W_{\so}(\check{g})$ we deduce the following formula: $W_{\su}(\check{g})=4W_{\so}(g)$. 

Dividing the equation \eqref{produit} by 4, we obtain the conclusion of the lemma.
\end{proof}

One can define now the extension of  $\Gamma_{\so}$ that we are looking for.
\begin{prop}
Let $\hat{\Gamma}$ be the set of pairs $(g,u)$ where $g\in\Gamma_{\so}$ and $u\in\R/\Z$ is such that $u=W_{\so}(g)\mod \frac{1}{4}\Z$. The operation $$(g,u)\cdot (h,v)=(gh,u+v-\frac{1}{2}\int_{\Sigma}\lg g^{-1}\d g \wedge \d hh^{-1}\ld)$$ induces a group structure on $\hat{\Gamma}$. Moreover this group acts on $\Omega^1_{\rm flat}\times S^1$ by the formula $$(a,\theta)^{(g,u)}=(a^g,\theta+\frac{1}{2}\int_{\Sigma}\lg g^{-1}ag\wedge g^{-1}\d g\ld -u).$$
\end{prop}
\begin{proof}
The only non obvious fact is that the product $(g,u)\cdot(h,v)$ is indeed an element of $\hat{\Gamma}$, but this is a direct consequence of the lemma \ref{lemme1}.
\end{proof}

\subsubsection{Properties of $\hat{\Gamma}$}
There is a natural map from $\Gamma_{\su}$ to $\hat{\Gamma}$ sending $g$ to the couple $(g,W_{\su}(g))$. The kernel of this map consists in constant functions with value $\pm 1$. By an abuse of notation, we will identify an element and its image.
We have the following proposition:

\begin{prop}\label{es}
The image of $\Gamma_{\su}$ in $\hat{\Gamma}$ is normal: let $\hat{H}$ denote the quotient. This group acts as a transformation group of $\boL$ and fits in the following exact sequence:
$$0\to\Z/4\to\hat{H}\to H^1(\Sigma,\Z/2)\to 0$$
\end{prop}
 
\begin{proof}
Let $h$ belong to $\Gamma_{\su}$ and the pair $(g,u)$ belong to $\hat{\Gamma}$. The conjugation $(g,u)\cdot h\cdot (g,u)^{-1}$ is a pair of the form $(ghg^{-1},V(g,u,h))$ where $V$ is a continuous function of $g,u$ and $h$. As the class of $h$ is 0 in $H^1(\Sigma,\Z_2)$, the same is true for $ghg^{-1}$, hence this function lifts to $\su$ and we still have the equality $V(g,u,h)=W_{\su}(ghg^{-1})\mod \frac{1}{4}\Z$. This equality is true mod $\Z$ when $h$ is the identity and in general by continuity as we can connect $h$ to the identity with a path. Finally we proved that $(g,u)\cdot h\cdot (g,u)^{-1}$ is in the image of $\Gamma_{\su}$.

We define a map from $\hat{\Gamma}$ to $H^1(\Sigma,\Z_2)$ by sending a pair $(g,u)$ to the homology class $g^*\alpha$ where $\alpha$ is the generator of $H^1(\so,\Z_2)$. In the case it is  0, $g$ lifts to $\su$ and $u$ has the form $W_{\su}(g)+x/4$ for $x\in\Z_4$. Elements of $\Gamma_{\su}$ are precisely those for which $x=0$. This explains the exact sequence of the proposition. 
\end{proof}

For any oriented curve $\gamma$ in $\Sigma$, we define an element $[\gamma]\in \hat{H}$ in the following way: suppose that $\Phi: S^1\times [0,1]\to\Sigma$ is an oriented embedding such that $\gamma=\Phi(S^1\times\{0\})$ and that the orientation of $\gamma$ coincides with the orientation of $S^1\times\{0\}$. We define an element $g_{\gamma}$ of $\Gamma_{\so}$ as being equal to $1$ in the complement of the image of $\Phi$ and satisfying the following equation;
$$\forall(t,s)\in S^1\times [0,1], \quad g_{\gamma} (\Phi(t,s))=\exp(i\pi\int_0^s\phi(x)\d x).$$
In this formula, $\phi$ is a smooth function with support in $[0,1]$ whose integral is equal to 1.

\begin{defn}\label{invol}
For any oriented curve $\gamma$ as above, we define an element $[\gamma]\in \hat{H}$ as the class of the pair $(g_{\gamma},0)$. The image of $[\gamma]$ in $H^1(\Sigma,\Z_2)$ in the exact sequence of the proposition \ref{es} is equal to $\gamma^{\#}$ the Poincaré dual of $\gamma$.
\end{defn}
This definition makes sense because $W_{\so}(g_{\gamma})=0 \mod \frac{1}{2}\Z$. This is a consequence of the following useful remark: 
any function $g$ with values in $\exp(i\R)\subset \su$ will satisfy $W_{\so}(g)=0$ as we can find a 3-manifold $M$ and an extension $\tilde{g}$ of $g$ on $M$ which take again values in $\exp(i\R)$. We obtain easily that $\tilde{g}^*\chi=0$ and the remark follows.

These elements are sufficient to understand the group structure of $\hat{H}$ as it is shown in the next proposition: we obtain the same Heisenberg group as in \cite{BHMV} and \cite{am}.

\begin{prop}
For any oriented curve $\gamma$ in $\Sigma$, the element $[\gamma]$ depends only on the class of $\gamma$ in $H_1(\Sigma,\Z)$. Denote by $\tau$ the class of the pair $(1,\frac{1}{4})\in \hat{\Gamma}$. Then $\hat{H}$ is generated by the elements $[\gamma]$ and $\tau$ with the following relations:
\begin{itemize}
\item[-]
$\tau^4=1$ and $[\gamma]^2=1$.
\item[-]
For all $\gamma$ and $\delta$ one has $[\gamma]\cdot [\delta]=\tau^{-\gamma\cdot\delta}[\epsilon]$ where $\epsilon$ represents the sum of the homology classes of $\gamma$ and $\delta$.
\end{itemize}
\end{prop}
\begin{proof}
If $g$ is an element of $\Gamma_{\so}$ with values in $\exp(i\R)$, we denote by ${\bf g}$ the function with values in $\R/\Z$ such that $g=\exp(i\pi{\bf g})$. Abusing notation, we will identify $g$ and ${\bf g}$: the formula for the product becomes:
$$({\bf g},u)\cdot ({\bf h},v)=({\bf g}+{\bf h},u+v-\frac{1}{4}\int_{\Sigma}\d{\bf g}\wedge\d{\bf h}).$$ 

Let $\gamma$ and $\delta$ be homologous curves. We remark that ${\bf g}_{\gamma}$ and ${\bf g}_{\delta}$ as functions from $\Sigma$ to $S^1$ represent the integral cohomology classes $\gamma^{\#}$ and $\delta^{\#}$, and hence they are homotopic. In other terms, the difference ${\bf h}={\bf g}_{\delta}-{\bf g}_{\gamma}$ lifts to a function from $\Sigma$ to $\R$. Then we have $({\bf h},0)\cdot ({\bf g}_{\gamma},0)=({\bf g}_{\delta},\frac{-1}{4}\int_{\Sigma}\d{\bf h}\wedge \d{\bf g}_{\gamma})$. 
But one has $\d{\bf h}\wedge \d{\bf g}_{\gamma}=\d\left({\bf h}\wedge \d{\bf g}_{\gamma}\right)$ as ${\bf h}$ has real values. By Stokes theorem, the integral vanishes. We conclude that the equality $[\gamma]=[\delta]$ is verified.

Now, let $\gamma$ and $\delta$ be any curves. By construction $\d{\bf g}_{\gamma}$ is Poincar\'e dual to $\gamma$ and the same is true for $\delta$. Hence, we have 
$\int_{\Sigma}\d{\bf g}_{\gamma}\wedge \d{\bf g}_{\delta}=\gamma\cdot\delta$. As ${\bf g}_{\gamma}+{\bf g}_{\delta}$ represents the sum of the cohomology classes $\gamma^{\#}$ and $\delta^{\#}$, this concludes the proposition. \end{proof}

\subsection{Hamiltonian flows and their holonomies}
\subsubsection{Normal forms of flat connections along curves}
In order to decribe the twist flows and compute their holonomies, we will use the following lemma:
\begin{lem}
Let $\Sigma$ be an oriented surface, $\Phi:S^1\times [0,1]\to\Sigma$ an oriented embedding and $a\in \Omega^1(\Sigma,\bog)$ a flat connection whose holonomy along $\Phi(S^1\times\{0\})$ differs from $\pm 1$. Let $t$ be the coordinate of $S^1$ viewed as $\R/\Z$.
Then there is a neighborhood $U$ of $a$ in $\Omega^1(\Sigma,\bog)$, a continuous function $g:U\to\Gamma_{\su}$ and a continuous function $\xi:U\to\bog$ such that for all $b\in U$ one has:
$$\Phi^* b^{g(b)}=\xi(b)\d t$$

\end{lem}
In other terms, up to gauge transformation, one can normalize a flat connection in the neighborhood of a curve such that in local coordinates, it looks like $\xi\d t$ where $t$ is a local $S^1$ parameter. Moreover, if $\exp(\xi)\ne \pm 1$, one can normalize in a continuous way all connections around the original one. The non-parametric lemma is part of the folklore, and the parametric version is an easy adaptation so we omit the proof.

\subsubsection{Description of the twist flows}

Let $\Sigma$ be a surface and $\gamma$ a curve in $\Sigma$ not homotopic to 0. We consider the gauge invariant function $F_{\gamma}:\Omega^1_{\rm flat}(\Sigma,\bog)\to [0,1]$ defined by $F_{\gamma}(a)=\frac{1}{\pi}\acos \re \hol_{\gamma} a$. We recall that the real part of a quaternion is equal to half the trace of the corresponding element in $\su$.
By twist flow, we mean the hamiltonian flow of the function $F_{\gamma}$. There is a problem in its definition as the moduli space $\boM(\Sigma,\su)$ is smooth only at irreducible representations. We will prove that the flow itself is well defined on $F_{\gamma}^{-1}(]0,1[)$ and that it is actually the hamiltonian flow of $F_{\gamma}$ at each regular point.

\begin{prop}
Let $\Sigma$ be a surface and $\Phi:S^1\times [0,1]\to\Sigma$ an oriented embedding around the curve $\gamma=\Phi(S^1\times\{0\})$. We denote by $t$ the coordinate of $S^1=\R/\Z$ and by $s$ the coordinate of $[0,1]$.

Let $a\in \Omega^1(\Sigma,\bog)$ be a flat connection normalized in the image of $\Phi$, that is such that there is $\xi\in\bog$ satisfying $\Phi^*a=\xi\d t$. We suppose that $\exp(\xi)\ne \pm 1$, or in other terms $||\xi||\notin \N \pi$ where $||\xi||^2=(\xi,\xi)$. Let $\phi$ be a smooth function with support in $[0,1]$ whose integral is equal to 1. 

Then for any $T\in \R$ we set $a^T=a+2\pi T\phi(s)\frac{\xi}{||\xi||}\d s.$

The projection of $a^T$ on $\boM(\Sigma,\su)$ defines a flow which depends only on the class of $a$ in $\boM(\Sigma,\su)$. This flow is the hamiltonian flow of $F_{\gamma}$ at every regular point of $\boM(\Sigma,\su)$.
\end{prop}

\begin{proof}
Let us prove first that this flow does not depend on the normalization. If $a$ can be normalized in two different ways, this means that there is a gauge transformation $g\in\Gamma_{\su}$ and $\xi,\eta$ such that $\exp(\xi)$ and $\exp(\eta)$ belong to the same non-central conjugacy class and the equation $g^* (\xi \d t)=\eta \d t$ is satisfied in the chart defined by $\Phi$. This equation implies that $g$ is a function of $t$ and satisfies the equation $\frac{\d g}{\d t}=g\eta-\xi g$.
The solution of this equation is $g(t)=\exp(-t\xi)g(0)\exp(t\eta)$ with the condition that $g(1)=g(0)$. A direct computation shows that this gauge transformation satisfies $g^* a^T=\eta \d t+2\pi T\phi(s)\frac{\eta}{||\eta||}\d s$ in the chart defined by $\Phi$. Hence, the flows are the same in $\boM(\Sigma,\su)$.

Suppose that $a^T$ represents a regular point of $\boM(\Sigma,\su)$. In order to show that the tangent vector of the flow is the symplectic gradient of $F_{\gamma}$ one needs to show the following identity:
$$d_{a^T}F(b)+\int_{\Sigma}\lg v_T,b\ld=0\text{ for all }b\in T_{a^T}\Omega^1_{\rm flat}.$$
where $v_T=\frac{d}{dT}a^T=2\pi \frac{\xi}{||\xi||}\phi(s) ds$.
In order to prove this formula, it is sufficient to check it for $b$ belonging to a family of flat connections which generates the tangent space of the class of $a$ in $\boM(\Sigma,\su)$. Thanks to the parametric version of the normalization lemma, one can suppose that $b$ is normalized, that is, it has the form $\Phi^* b=\eta \d t$.
We have then $\pi F_{\gamma}(a^T+ub)=\acos\re \hol_{\gamma} (\xi dt+T\phi(s)\xi ds+u\eta dt)=\acos\re\exp(\xi+u\eta)=||\xi+u\eta|| \mod 2\pi$. So we have $d_{a^T}F(b)=\frac{(\xi,\eta)}{\pi ||\xi||}$.

Moreover one has $\int_{\Sigma}\lg v_T,b\ld=\int_{\Sigma}\lg 2\pi\phi(s)\frac{\xi}{||\xi||} ds, \eta dt\ld=-2\pi\lg\frac{\xi}{||\xi||},\eta\ld=-\frac{(\xi,\eta)}{\pi ||\xi||}$. These terms are opposite, and the proposition is proved.
\end{proof}

Using this expression of the twist flows, it is easy to see that the hamiltonian flow of $F_{\gamma}$ is 1-periodic. In the notation of the proposition, starting from a normalized connection $a$ i.e. such that $\Phi^*a=\xi \d t$, one obtains after a time 1 the connection $a+2\pi\phi(s)\d s \frac{\xi}{||\xi||}$. This connection is gauge equivalent to the first one: by setting  $g(t,s)=\exp(2i\pi\int_0^s\phi(x)\d x)$ in the chart given by $\Phi$ and extending $g$ by 1 outside the cylinder one has  $a^1=g^*a$. 

One can also interpret what happens at time 1/2: using the same notation, one has $a^{1/2}=a+\pi\phi(s)\d s \frac{\xi}{||\xi||}$ which is gauge equivalent to $a$ with gauge element in $\Gamma_{\so}$. The element $g$ satisfying $a^{1/2}=g^*a$ is precisely the one we used to define $[\gamma]$ in the definition \ref{invol}. By the discussion in the beginning of the section \ref{heisenberg}, one sees that the hamiltonian flow of $F_{\gamma}$ at time $1/2$ sends a representation $\rho$ to the product $\rho \gamma^{\#}$ where $\gamma^{\#}$ is the Poincar\'e dual of $\gamma$ in $H^1(\Sigma,\Z_2)$, viewed as a representation with values in $\{\pm 1\}$. We recover the fact that the twist flows along separating curves have period $1/2$ as these curves are homologous to 0.

\subsubsection{Holonomy of the twist flows}

All ingredients are present to compute the holonomy of the prequantum bundle along orbits of the twist flows. We will also compute the holonomy of half orbits using the action of the Heisenberg group. The results are summed up in the following proposition:

\begin{prop}
Let $\Sigma$ be a surface, and $\boL$ be the prequantum circle bundle over $\boM(\Sigma,\su)$ defined in the section \ref{generalite}. We will write the action of $S^1$ on $\boL$ on the right and the action of $\hat{H}$ on the left.
Let $\gamma$ be a non trivial curve and $\rho$ a representation in $\boM(\Sigma,\su)$. Suppose that $\rho(\gamma)\ne \pm 1$. Then 
\begin{itemize}
\item[-]
The parallel transport of an element $\alpha$ of $\boL$ over $\rho$ along the hamiltonian flow of $F_{\gamma}$ after a time 1 is equal to $\alpha-F_{\gamma}(\rho)$ and lies over $\rho$.
\item[-]
The parallel transport of an element $\alpha$ of $\boL$ over $\rho$ along the hamiltonian flow of $F_{\gamma}$ after a time 1/2 is equal to $[\gamma]\alpha-\frac{1}{2}F_{\gamma}(\rho)$ and lies over $\gamma^{\#}\rho$.
\end{itemize}
\end{prop}
\begin{proof}
As usual, we find a positive embedding $\Phi:S^1\times [0,1]\to \Sigma$ and suppose up to gauge transformation that there exists $\xi\in\bog$ such that one has $\Phi^* a=\xi\d t$. Hence, the twist flow may be described by the following 1-parameter family of flat connections : $a^T=a+2\pi T \phi(s)\d s \frac{\xi}{||\xi||}$. The gauge transformations which allow to compare $a,a^{1/2}$ and $a^1$ are given by the expressions $g_{\kappa}(t,s)=\exp(\kappa \pi \int_0^s\phi(x)\d x\frac{\xi}{||\xi||})$ where $\kappa=1$ if $T=1/2$ and $\kappa=2$ if $T=1$.

To compute the holonomy of $\boL$ along that path, we work on the product $\Omega^1(\Sigma,\bog)\times S^1$ and take $\theta$ such that $\alpha$ is represesented by $(a,\theta)$. We need to find a function $\theta^T$ such that $\theta^0=\theta$ and $\lambda(\frac{\d}{\d t}(a^T,\theta^T))=0$ for all $T$ where $\lambda$ is the connection 1-form defined in section \ref{generalite}. After simplification, one gets 
 $\frac{||\xi||}{2\pi}+\frac{\d \theta^T}{\d T}=0$ which gives $\theta^T=-\frac{T}{2\pi}||\xi||$. 
 We can compare these elements for $T\in \{0,1/2,1\}$ by using the action of the extended gauge group on  $\Omega^1(\Sigma,\bog)\times S^1$. Let us do it precisely for $T=1/2$.
 The element $(g_1,0)\in\hat{\Gamma}$ acts on $a$ by the formula $(a,\theta)^{(g_1,0)}=(a^{1/2},\theta+\frac{1}{2}\int_{\Sigma}\lg g_1^{-1}a g_1\wedge g_1^{-1}\d g_1\ld -0)=(a^{1/2},\theta+\frac{||\xi||}{4\pi})$.
 
 We deduce from this formula that $(a^{1/2},\theta^{1/2})=(a,\theta)^{(g_1,0)}-\frac{||\xi||}{2\pi}$. When projecting onto $\boL$ one finds that the class of the element $(a,\theta)$ has moved to $[\gamma](a,\theta)-\frac{1}{2}F_{\gamma}(\rho)$, which was claimed in the proposition. This implies the case of time 1 as one has $[\gamma]^2=1$.
\end{proof}

\subsection{Identifications}

We are now ready to make the link between the two parts and explain the relation between the algebra $K(\Sigma,-i)$ and the algebra generated by parallel transport along twist flows during time $\frac{1}{2}$.
Recall that we defined $\Phi_{\gamma}^t$ to be the hamiltonian flow of $F_{\gamma}$ extended to $\boL$. It extends also to the hermitian line bundle $\boL^{\C}$ associated  to $\boL$.
Again, this flow is well-defined on the preimage of $]0,1[$ by $F_{\gamma}$. From it, one can define a linear operator $\Psi_{\gamma}^t$ on $L^2(\boM(\Sigma,\su),\boL^{\C})$ thanks to the formula $(\Psi_{\gamma}^t s)(\rho)=\Phi_{\gamma}^t s(\Phi_{\gamma}^{-t} \rho)$. 
By defining $\boO_{\gamma}=\Psi_{\gamma}^{1/2}+\Psi_{\gamma}^{-1/2}$ one obtains now easily the following proposition:

\begin{prop}
The map from $K(\Sigma,-i)$ to $\en L^2(\boM(\Sigma,\su),\boL)$ sending $\gamma$ to $\boO_{\gamma}$ is an  injective morphism of algebras.
\end{prop}
\begin{proof}
We will use the theorem \ref{combi} and replace $K(\Sigma,-i)$ with its isomorphic algebra $\boA(\Sigma)$. We recall that $\boA(\Sigma)$ is a sub-algebra of the tensor product $K(\Sigma,-1)\otimes \boA$ and remark that $\boA$ is equal to $\C\otimes_{\Z_4}\hat{H}$, where $\Z_4$ acts on $\C$ by $x.z=i^xz$.

The first observation is that $K(\Sigma,-1)$ is isomorphic to an algebra of functions on $\boM(\Sigma,\su)$ where a simple curve $\gamma$ gives rise to the function $f_{\gamma}(\rho)=-\tr \rho(\gamma)$. Hence there is an algebra homomorphism from $K(\Sigma,-1)$ to $\en L^2(\boM(\Sigma,\su),\boL^{\C})$ given by multiplication operators. On the other hand, the group $\hat{H}$ acts on $\boL$, and this action lies over the action of $H^1(\Sigma,\Z_2)$ on $\boM(\Sigma,\su)$ . Hence, the algebra of the group acts on sections of $\boL$ via the formula $([\gamma].s)(\rho)=[\gamma]s(\gamma^{\#}\rho)$. We recall that the $\Z_4$ central subgroup of $\hat{H}$ acts by rotation of angle $\pi/2$ on the fibres. By identifying the action of the central subgroup with the multiplication with $i$, one obtains finally an action of $\boA=\C\otimes_{\Z_4}\hat{H}$.
Joining the two actions, one obtains a natural morphism $\psi$ from $\boA(\Sigma)$ to $\en L^2(\boM(\Sigma,\su),\boL^{\C})$.

One computes 
\begin{eqnarray*}
(\boO_{\gamma}s)(\rho)&=&
\Phi_{\gamma}^{1/2} s(\Phi_{\gamma}^{-1/2} \rho)+\Phi_{\gamma}^{-1/2} s(\Phi_{\gamma}^{1/2} \rho)\\
&=&\exp(-i\pi F_{\gamma}(\rho))[\gamma]s(\gamma^{\#}\rho)+\exp(i\pi F_{\gamma}(\rho))[\gamma]s(\gamma^{\#}\rho)\\
&=&2\cos(\pi F_{\gamma}(\rho))[\gamma] s(\gamma^{\#}\rho)=\left(\tr \rho(\gamma)\right) [\gamma] s(\gamma^{\#}\rho)
\end{eqnarray*}
We interpret this formula by identifying $\boO_{\gamma}$ and $\psi((-1)^{n(\gamma)}\gamma\otimes [\gamma])$ which proves the proposition. The injectivity is clear.
\end{proof}
To conclude this section, we recall the following triangle of algebra identifications:

$$\xymatrix{
K(\Sigma,-i)\ar@/_1.5pc/[rr]_-{\simeq}^-{\phi}\ar[r]^-{\boO}&\boA(\Sigma)\ar[r]^-{\psi}
&\en L^2(\boM(\Sigma,\su),\boL^{\C})
}$$

\section{Traces}

Let $\frac{a}{b}$ be a rational number in irreducible form with $b>0$, $\zeta$ be a rational number and $p:\N^*\to 2\N^*$ be an increasing sequence such that for all $n>0$, the fraction $\theta_n=\frac{a}{b}+\frac{\zeta}{p_n}$ has $p_n$ as lowest denominator. We will call {\it admissible} such a sequence $\theta=(\theta_n)$.
We consider the topological quantum field theory (TQFT) $(V_n,Z_n)$ defined in \cite{BHMV} with root $A=\exp(i\pi\theta_n)$ which has order $2p_n$. It corresponds to the $\su$ WZW-model with level equal to $p_n/2-2$.

In \cite{mn}, it is shown that for $\theta_n=1+1/2n$ (i.e. for roots converging to -1), the representation of the mapping class group of $\Sigma$ on $\en V_n(\Sigma)$ converges in the Fell topology to the representation on the subspace of $L^2(\boM(\Sigma,\su))$ generated by trace functions. We will generalize this result in the following sense:
\begin{thm}
Let $\theta=(\theta_n)$ be an admissible sequence with limit $a/b$. We define a sesquilinear form on $K(\Sigma,e^{i\pi a/b})$ by the formula $\langle x,y\rangle=\langle xy\rangle_{\theta}$ where the product is understood in the algebra $K(\Sigma,e^{i\pi a/b})$ and the linear form $\langle\cdot\rangle_{\theta}$ is given in the proposition \ref{trace}.
Then, the sequence of representations of the mapping class group on $\en V_n(\Sigma)$ converges in the Fell topology to $K(\Sigma,e^{i\pi a/b})$.
\end{thm}
 Notice that we abuse the notion of Fell topology as the representations are not irreducible neither hermitian.
Nevertheless, for $\theta_n=-1/2+1/4n$ (i.e. for roots converging to $-i$), the limit representation is hermitian and is understood as the algebra of parallel transport operators described in the last section with a natural trace interpreted in the section \ref{geom}.
The proof of this theorem occupies the next section.

\subsection{Computation of traces in TQFT}
Let $\Sigma$ be a surface and $\gamma$ be a multicurve adapted to a pants decomposition. This pants  decomposition is parametrized by a pair $(H,G)$ where $H$ is a handlebody bounding $\Sigma$ and $G$ is a banded trivalent graph embedded in $H$ such that $H$ retracts on it. For any edge $e$ of $\Gamma$, we denote by $m_e$ the multiplicity of the curve dual to that edge in $\gamma$. The formula of the lemma 3.4 in \cite{mn} which is an easy computation using the basis provided in the theorem 4.11 of \cite{bhmv} may be adapted to our settings. One obtains the following formula:

\begin{equation}\label{somme}
\tr Z_n(\Sigma\times [0,1],\gamma)=\sum\limits_{\sigma} \prod\limits_e \left(-2\cos(2\pi\theta_n(\sigma_e+1))\right)^{m_e}
\end{equation}

In this formula, $e$ runs over edges of $G$ and $\sigma$ is a function from edges of $G$ to the set $\{0,1,\ldots,p_n/2-2\}$ which satisfies the following conditions for any triple of edges $(i,j,k)$ incident to the same vertex:
\begin{itemize}
\item[-] $\sigma_i+\sigma_j+\sigma_k$ is even and less or equal to $p_n-4$.
\item[-] $\sigma_i\le \sigma_j+\sigma_k$ and permutations of the same inequality.
\end{itemize}
Let $\Lambda$ be the group of integer valued functions on edges of $G$ such that the sum of the values of any three incident edges is even. 
In the sum \eqref{somme}, the map $\sigma$ runs over a subset of $\Lambda$. Let $B$ be the least common multiple of $b$ and 2. Then, one has $B\Z^{E(G)}\subset \Lambda$. We decompose the sum into classes modulo $B\Z^{E(G)}$, denoting by $\Lambda_B$ the quotient $\Lambda/B\Z^{E(G)}$
The sum may be written in the following way:
$\sum\limits_{\mu\in\Lambda_B}\sum\limits_{\ba{\sigma}=\mu} \prod\limits_e \left(-2\cos(2\pi\theta_n(\sigma_e+1))\right)^{m_e}$.

We recognize the second sum as a Riemannian one: more precisely, let $U_G$ be the set of functions $\tau$ from $E(G)$ to $[0,1]$ such that for any edges $(i,j,k)$ incident to the same vertex one has $\tau_i+\tau_j+\tau_k\le 2$, $\tau_i\le \tau_j+\tau_k$ and their cyclic permutations. This set is endowed with the Lebesgue measure $\d\tau$.

For any $\mu\in\Lambda_B$, let $F_{\mu}(\tau)=\prod\limits_e  (-2\cos(2\pi \frac{a}{b}(\mu_e+1)+\pi\zeta\tau_e))^{m_e}$ and $d_G$ be the number of edges of $G$. 
Then one has 
\begin{prop}\label{trace}
$$\lim\limits_{n\to\infty}(\frac{2}{p_n})^{d_G}\tr Z_n(\Sigma\times [0,1],\gamma)=\frac{1}{B^{d_G}}\sum\limits_{\mu\in\Lambda_B} \int_{U_G}F_\mu(\tau)\d\tau.$$
We denote this limit as $\langle\gamma\rangle_{\theta}$ as it depends only on $\gamma$ and on the admissible sequence $\theta=(\theta_n)$.
\end{prop}

Let us explicit the proposition when $a/b=-1/2$, $\zeta=1$ and $p_n=4n$.

\begin{eqnarray*}
\lim\limits_{n\to\infty}(\frac{2}{p_n})^{d(G)}\tr Z_n(\Sigma\times[0,1],\gamma)&=\frac{1}{2^{d_G}}
\sum\limits_{\mu}\int\limits_{U}\prod\limits_e(-2\cos(\pi \mu_e+\pi+\pi\tau_e))^{m_e}\d\tau \\
&=\frac{1}{2^{d_G}}(\sum\limits_{\mu}(-1)^{\sum_e \mu_e m_e})\int\limits_{U}\prod\limits_e(2\cos(\pi\tau_e))^{m_e}\d\tau
\end{eqnarray*}
In this formula, $\mu\in\Lambda_2\simeq H_1(G,\Z_2)$. The tuple $(m_e)$ may be viewed as an element of $H^1(G,\Z_2)$. Hence, the sum $\sum\limits_{\mu}(-1)^{\sum_e \mu_e m_e}$ vanishes if the class of $(m_e)$ is non zero and is equal to $2^g$ in the other case, where $g$ is the dimension of $H^1(G,\Z_2)$.

We remark that the class of $(m_e)$ vanishes in $H^1(G,\Z_2)$ if and only if the class of $\gamma$ noted $\gamma^{\#}$ vanishes in $H^1(\Sigma,\Z_2)$. Moreover, the integral may be interpreted as in \cite{mn} as the integral of the function $(-1)^{n(\gamma)}f_{\gamma}$ on $\boM(\Sigma,\su)$ with respect to the Liouville form $\d \lambda$.

Finally, we can write the following formula:
\begin{equation}\label{tracei}
\lim\limits_{n\to\infty}(\frac{2}{p_n})^{d(G)}\tr Z_n(\Sigma\times[0,1],\gamma)=\delta_{\gamma^{\#}=0}\int_{\boM(\Sigma,\su)}(-1)^{n(\gamma)}f_{\gamma}d\lambda.
\end{equation}
\begin{rem}
Such a formula may follow from the tensor product formula of TQFT for even $p_n$ (see thm. 1.6 in \cite{BHMV}) but we did not go on in that direction.
\end{rem}

\subsection{Geometric interpretation of the trace}\label{geom}

By construction of the TQFT of \cite{BHMV}, the limit involved in the proposition \ref{trace} gives a trace on $K(\Sigma,\pm\exp(\frac{i\pi a}{b}))$: let $\tr$ be the linear form defined by 
$$\tr(\gamma)=\lim\limits_{n\to\infty}(\frac{2}{p_n})^{d(G)}\tr Z_n(\Sigma\times[0,1],\gamma).$$
 We recall that  this linear form depends on $\zeta$ and that it satisfies $\tr(\gamma\delta)=\tr(\delta\gamma)$ for any multicurves $\gamma$ and $\delta$. 

This trace is well understood in the case of $K(\Sigma,-1)$: one has $\tr(\gamma)=\int_{\boM(\Sigma,\su)}f_{\gamma}\d\lambda$. We propose to explain the geometry of the formula \eqref{tracei}.

Let $\tr:\boA\to\C$ be the linear form defined by $\tr([\gamma])=\delta_{\gamma^{\#}=0}$ where the equality $\gamma^{\#}=0$ is understood in $H^1(\Sigma,\Z_2)$. Viewed in $\boA(\Sigma)$, the trace defined in the formula \eqref{tracei} is the tensor product of the traces we have just recalled on $K(\Sigma,-1)$ and $\boA$. This is a clear geometric construction of the trace but we find even more geometric to think of it in terms of parallel transport operators on $L^2(\boM(\Sigma,\su),\boL^{\C})$. 

All operators involved in this article have the form $(\boO s)(\rho)=K(\rho)s(f(\rho))$ where $f$ is a transformation of $\boM(\Sigma,\su)$ which is either the identity or with a fixed point set of codimension grater than 1, and $K$ is a section of the bundle $\hom({\it f}^* \boL^{\C},\boL^{\C})$. Mimicking the definition of the trace of an operator with kernel, it is natural to define $\tr \boO=0$ if $f$ has a measure zero set of fixed points and $\tr \boO= \int_{\boM(\Sigma,\su)}(\tr K) \d \lambda$ if $f$ is the identity. The trace we obtain in this way by interpreting  $K(\Sigma,-i)$ as an algebra of operators is exactly the same trace we have defined previously.

We would like to end with some questions:
\begin{itemize}
\item[-]
If we associate to a curve $\gamma$ the operator $\boO_{\gamma}=\Psi^{a/b}+\Psi^{-a/b}$ on $L^2(\boM,\boL^{\C})$, we obtain formulas for the geometric traces which look very similar to the traces on $K(\Sigma,\exp(\frac{i\pi a}{b}))$ obtained by taking limits of TQFT traces. 
Unfortunately, we could not identify exactly the formulas and we did not find a representation of the skein algebra in $L^2(\boM,\boL^{\C})$ which would give sense to these equalities. Still, this coincidence is very fascinating, is there any convincing explanation?
\item[-]
At first, we did not expect the trace on $K(\Sigma,\pm i)$ to be positive, as it appears as a limit of non hermitian TQFT but it is actually. Is it true that all limit traces on $K(\Sigma,\pm\exp(i\pi\frac{a}{b}))$ are positive?
\end{itemize}

\textsc{UPMC Univ Paris 06, UMR 7586, Institut de Mathématiques de Jussieu, F-75005, Paris, France}\\
{\it E-mail adress:} {\tt marche@math.jussieu.fr}

\end{document}